\documentclass[a4paper,12pt]{article}
\usepackage{color}
\usepackage[english]{babel}
\usepackage{amsfonts}
\usepackage{amsmath}
\usepackage{amsthm}
\usepackage{amssymb}
\usepackage{bbm}
\usepackage{mathrsfs}
\usepackage{esint}
\usepackage{faktor}
\usepackage[utf8]{inputenc}
\usepackage{afterpage}
\usepackage{graphicx}

\newcommand{\R}{\mathbb{R}}

\renewcommand{\epsilon}{\varepsilon}
\renewcommand{\phi}{\varphi}

\newtheorem{lemma}{Lemma}[section]
\newtheorem{thm}[lemma]{Theorem}

\newtheorem{cor}[lemma]{Corollary}

\theoremstyle{definition}

\newtheorem{rmk}[lemma]{Remark}

\numberwithin{equation}{section}

\DeclareMathOperator*{\essinf}{ess \, inf}

\begin{document}
\title{\textbf{Some recent results on singular $p$-Laplacian equations}}
\author{
\bf Umberto Guarnotta\thanks{Corresponding author.}, Roberto Livrea\\
\small{Dipartimento di Matematica e Informatica, Universit\`a di Palermo,}\\
\small{Via Archirafi 34, 90123 Palermo, Italy}\\
\small{\it E-mail: umberto.guarnotta@unipa.it, roberto.livrea@unipa.it}\\
\mbox{}\\
\bf Salvatore A. Marano\\
\small{Dipartimento di Matematica e Informatica, Universit\`a di Catania,}\\
\small{Viale A. Doria 6, 95125 Catania, Italy}\\
\small{\it E-mail: marano@dmi.unict.it}\\
\mbox{}\\
}
\date{}
\maketitle
\begin{abstract}
A short account of some recent existence, multiplicity, and uniqueness results for singular $p$-Laplacian problems either in bounded domains or in the whole space is performed, with a special attention to the case of convective reactions. An extensive bibliography is also provided.
\end{abstract}
\vspace{2ex}
\vspace{2ex}


\section{Introduction}
When studying quasi-linear elliptic systems in the whole space and with singular, possibly convective, reactions, a natural preliminary step is looking for the previous literature on equations of the same type, which we have done in the latest years.
 
At first, this obviously led us to investigate \textit{singular $p$-Laplacian Dirichlet problems} as
\begin{equation}\label{sample1}
\left\{
\begin{alignedat}{2}
-\Delta_p u & =h(x,u,\nabla u)\;\; &&\mbox{in}\;\;\Omega,\\
u & >0 &&\mbox{in}\;\;\Omega,\\
u &= 0 &&\mbox{on}\;\;\partial\Omega,
\end{alignedat}
\right.
\end{equation}
where $1<p<\infty$, the symbol $\Delta_p$ denotes the $p$-Laplace operator, namely 
\begin{equation*}
\Delta_p u:={\rm div}(|\nabla u|^{p-2}\nabla u), 
\end{equation*}
$\Omega$ is a bounded domain in $\R^N$, $N\geq 3$, with smooth boundary $\partial\Omega$, while $h\in C^0(\Omega\times\R^+\times\R^N)$ satisfies
$$\lim_{t\to 0^+}h(x,t,\xi) =\infty.$$
If $p=2$ then various special (chiefly non-convective) cases of \eqref{sample1} have been thoroughly studied (see Subsection \ref{bound.1}). Both surveys \cite{HMV,HM,Ra} and a monograph \cite{GR}, besides many proceeding papers, are already available. The main purpose of Section \ref{bound} below is to provide a short account on some recent existence, multiplicity or uniqueness results for $p\neq 2$ and the relevant technical approaches. Let us also point out \cite{Sao,PVV,PRR}. The work \cite{Sao} treats a singular $p(x)$-Laplacian Robin problem, while \cite{PVV,PRR} are devoted to singular $(p,q)$-Laplacian equations with Neumann and Robin boundary conditions, respectively; cf. \cite{PRR20} too.

Section \ref{whole} aims at performing the same as regards \textit{singular $p$-Laplacian problems in the whole space}. So, it deals with situations like
\begin{equation}\label{sample2}
\left\{
\begin{alignedat}{2}
-\Delta_p u & =h(x,u,\nabla u)\;\; &&\mbox{in}\;\;\R^N,\\
u & >0 &&\mbox{in}\;\;\R^N,\\
u & (x) \to 0 &&\mbox{as}\;|x|\to\infty.
\end{alignedat}
\right. 
\end{equation}
To the best of our knowledge, except \cite{GR}, even when $p=2$ and $h$ does not depend on $\nabla u$, there are no surveys concerning \eqref{sample2}. Hence, this probably represents the first contribution.

Both sections are divided into four parts. The first is a historical sketch of the case $p=2$. The next two treat existence, multiplicity, and uniqueness in the non-convective case. The fourth is devoted to singular problems with convection. 
Since the literature on \eqref{sample1}--\eqref{sample2} is by now very wide and our knowledge is limited, significant works may have been not mentioned here, something of which we apologize in advance. Moreover, for the sake of brevity, we did not treat singular parabolic boundary-value problems and refer the reader to \cite{HM,DeBG,OlP2,GKS,CG}.
\section{Basic notation}
Let $X(\Omega)$ be a real-valued function space on a nonempty measurable set $\Omega\subseteq\R^N$. If $u_1,u_2\in X(\Omega)$ and $u_1(x)<u_2(x)$ a.e. in $\Omega$ then we simply write $u_1<u_2$. The meaning of $u_1\leq u_2$, etc. is analogous. Put 
\begin{equation*}
X(\Omega)_+:=\left\{u\in X(\Omega): u\geq 0\right\}.   
\end{equation*}
The symbol $u\in X_{loc}(\Omega)$ means that $u:\Omega\to\R$ and $u\lfloor_K\in X(K)$ for all nonempty compact subset $K$ of $\Omega$. Given $1<r<N$, define
$$r':=\frac{r}{r-1}\, , \quad r^*:=\frac{Nr}{N-r}\, .$$
%
%
Let us next recall the notion and some relevant properties of the so-called Beppo Levi space
$\mathcal{D}^{1,r}_0(\R^N) $, addressing the reader to \cite[Chapter II]{Ga} for a complete treatment. Set
\begin{equation*}
\mathcal{D}^{1,r} := \left\{ z \in L^1_{\rm loc}(\R^N): |\nabla z| \in L^r(\R^N) \right\}
\end{equation*}
and denote by $ \mathcal{R} $ the equivalence relation that identifies two elements in $ \mathcal{D}^{1,r} $ whose difference is a constant. The quotient set $\mathcal{\dot D}^{1,r}$, endowed with the norm
$$\|u\|_{1,r} :=\left(\int_{\R^N}|\nabla u(x)|^r{\rm d}x\right)^{1/r},$$
turns out complete. Write $\mathcal{D}^{1,r}_0(\R^N)$ for the subspace of $ \mathcal{\dot D}^{1,r} $ defined as the closure of $ C^\infty_0(\R^N) $ under $ \|\cdot\|_{1,r} $, namely
$$\mathcal{D}^{1,r}_0(\R^N) := \overline{C^\infty_0(\R^N)}^{\|\cdot\|_{1,r}}.$$
$\mathcal{D}^{1,r}_0(\R^N)$, usually called Beppo Levi space,  is reflexive and continuously embeds in $L^{r^*}(\R^N)$, i.e., 
\begin{equation}\label{embedding}
\mathcal{D}^{1,r}_0(\R^N) \hookrightarrow L^{r^*}(\R^N).
\end{equation}
Consequently, if $u\in\mathcal{D}^{1,r}_0(\R^N) $ then $u$ vanishes at infinity, meaning that the set $\{x\in\R^N:|u(x)|\geq\epsilon\}$ has finite measure for any $ \epsilon > 0 $. 
%
%
%
%
%
%
\section{Problems in bounded domains}\label{bound}
\subsection{The case $p=2$}\label{bound.1}
Let $\Omega$ be a bounded domain in $\R^N$, $N\geq 3$, with smooth boundary $\partial\Omega$, let $a:\Omega\to\R^+_0$ be nontrivial measurable, and let $\gamma>0$. The simplest singular elliptic Dirichlet problem writes as
\begin{equation}\label{prob0}
\left\{
\begin{alignedat}{2}
-\Delta u & =a(x)u^{-\gamma}\;\; &&\mbox{in}\;\;\Omega,\\
u & >0 &&\mbox{in}\;\;\Omega,\\
u &= 0 &&\mbox{on}\;\;\partial\Omega.
\end{alignedat}
\right.
\end{equation}
Since the pioneering papers \cite{S,CRT,CoP,LM,CLMK}, a wealth of existence, uniqueness or multiplicity, and regularity results concerning \eqref{prob0} have been published. We refer the reader to the monograph \cite{GR} as well as the surveys \cite{HMV,HM} for an exhaustive account. Roughly speaking, four basic questions can be identified:
\begin{itemize}
\item find the right conditions on the datum $a$. Usually, $a\in L^q(\Omega)$ with $q\geq1$ is enough for existence. However, starting from the works \cite{OP,OlP}, the case when $a$ is a bounded Radon measure took interest.
\item consider non-monotone singular terms. This is a difficult task, mainly when we want to guarantee uniqueness of solutions.
\item insert convective terms on the right-hand side. For equations driven by the Laplacian, good references are \cite[Section 9]{GR} and \cite{AL}. Otherwise, cf. \cite{Bo,Aalt,GMM,PaZ}.
\item substitute the Laplacian with more general elliptic operators. Obviously, a first attempt might be considering equations driven by the $p$-Laplacian, and this section aims to provide a short account of the nowadays literature. However, further possibly non-homogeneous operators have been considered; see, e.g., \cite{S,CRT,Go,Chab,Cui,BO1,OlP,GKSzero,PW2}.
\end{itemize}
Incidentally, we recall that \eqref{prob0} stems from important applied questions, as the study of heat conduction in electrically conducting materials \cite{FM}, chemical heterogeneous catalysts \cite{Pe}, and non-Newtonian fluids \cite{AM}.
\subsection{Existence and multiplicity}\label{S2.1}
Consider the model problem
\begin{equation}\label{prob1}
\left\{
\begin{alignedat}{2}
-\Delta_p u & =a(x)u^{-\gamma}+\lambda f(x,u) \; \; &&\mbox{in} \; \;\Omega, \\
u & >0 &&\mbox{in} \; \;\Omega, \\
u &= 0 &&\mbox{on} \; \;\partial\Omega,
\end{alignedat}
\right.
\end{equation}
where $a:\Omega\to\R^+_0$ denotes a nonzero measurable function, $\gamma,\lambda>0$, while $f:\Omega\times\R^+_0\to\R$ satisfies Carath\'eodory's conditions. Let us stress that, here, \textit{the parameter $\lambda$ multiplies the non-singular term}.

In 2006, Perera and Silva investigated \eqref{prob1} under the assumptions below, where $f$ is allowed to change sign.
\begin{itemize}
\item[$({\rm a}_1)$] There exist $\varphi_0\in C^1_0(\overline{\Omega} )_+$ and $\hat q>N$ such that $a\varphi_0^{-\gamma}\in L^{\hat q}(\Omega)$.
\item[$({\rm a}_2)$] With appropriate $\delta,c_1>0$ one has
$$f(x,t)\geq-c_1 a(x)\;\;\mbox{in}\;\;\Omega\times [0,\delta].$$
\item[$({\rm a}_3)$] To every $M>0$ there correspond $h\in L^1(\Omega)$ and $c_2>0$ such that
$$-h(x)\leq f(x,t)\leq c_2\quad\forall\, (x,t)\in\Omega\times [0,M].$$
\item[$({\rm a}_4)$] With appropriate $q\in ]1,p^*[$ and $c_3>0$ one has
$$f(x,t)\leq c_3(t^{q-1}+1)\;\;\mbox{in}\;\;\Omega\times\R^+_0.$$
\item[$({\rm a}_5)$] There are $t_0>0$ and $\mu>p$ such that
$$0<\mu\int_0^tf(x,\tau){\rm d}\tau\leq tf(x,t)\quad\forall\, (x,t)\in\Omega\times[t_0,+\infty[\, .$$
\end{itemize}
They seek distributional solutions to \eqref{prob1}, i.e., functions $u\in W^{1,p}_0(\Omega)$ such that $u>0$ and
\begin{equation*}
\int_\Omega |\nabla u|^{p-2}\nabla u\cdot\nabla\varphi{\rm d}x=\int_\Omega a u^{-\gamma}\varphi{\rm d}x+\int_\Omega f(\cdot, u)\varphi{\rm d}x\quad\forall\,\varphi\in C^\infty_0(\Omega).
\end{equation*}
\begin{thm}[\cite{PS}, Theorems 1.1--1.2]\label{thm1}
Let $({\rm a}_1)$--$({\rm a}_3)$ be satisfied. Then Problem \eqref{prob1} admits a distributional solution for every $\lambda>0$ small. If, in addition, $({\rm a}_4)$--$({\rm a}_5)$ hold true then a further distributional solution exists by decreasing $\lambda$ when necessary.
\end{thm}
Proofs employ perturbation arguments and variational methods, previously introduced in \cite{PZ}. An immediate but hopefully useful consequence of Theorem \ref{thm1} is the next 
\begin{cor}
Let $({\rm a}_1)$ be fulfilled. Suppose $f$ does not depend on $x$ and, moreover, $f(t)\geq 0$ in a neighborhood of zero once $\essinf_{\Omega}a=0$. Then, for every $\lambda>0$ sufficiently small, the problem
\begin{equation}\label{prob2}
-\Delta_p u=a(x)u^{-\gamma}+\lambda f(u) \; \;\mbox{in} \; \;\Omega,\quad
u>0\;\;\mbox{in} \; \;\Omega,\quad
u= 0\;\;\mbox{on} \; \;\partial\Omega
\end{equation}
possesses a distributional solution.
\end{cor}
Further results concerning \eqref{prob2} can be found in Aranda-Godoy \cite{AG}, where a continuous non-increasing function $g(u)$ takes the place of $u^{-\gamma}$ and, from a technical point of view, fixed point theorems for nonlinear eigenvalue problems are exploited. 

The case $\lambda=0$ in \eqref{prob2} was well investigated by Canino, Sciunzi, and Trombetta \cite{CST}, with a special attention to uniqueness (see the next section). Here, given $u\in W^{1,p}_{loc}(\Omega)$,
$$u=0\;\;\mbox{on}\;\;\partial\Omega\stackrel{\rm def}{\iff} u\geq 0\;\;\mbox{and}\;\; (u-\epsilon)^+\in W^{1,p}_0(\Omega)\;\;\forall\,\epsilon>0.$$
\begin{thm}[\cite{CST}, Theorem 1.3]\label{CSTthm1.3}
Let $\lambda=0$. If $\gamma\geq 1$ and $a\in L^1(\Omega)$ then \eqref{prob2} admits a distributional solution $u\in W^{1,p}_{loc}(\Omega)$ such that $\essinf_{K}u>0$ for any compact set $K\subseteq\Omega$. Moreover, $u^{1+(\gamma-1)/p}\in W^{1,p}_0(\Omega)$. If $0<\gamma<1$ then \eqref{prob2} has a solution $u\in W^{1.p}_0(\Omega)$ in each of the following cases:
\begin{itemize}
\item $1<p<N$ and $a\in L^m(\Omega)$, with $m:=\left(\frac{p^*}{1-\gamma}\right)'$.
\item $p=N$ and $a\in L^m(\Omega)$ for some $m>1$.
\item $p>N$ and $a\in L^1(\Omega)$. 
\end{itemize}
\end{thm}
The proof of this result relies on a technique previously introduced in \cite{BO1} for the semi-linear case. It employs truncation and regularization arguments. The work \cite{GCS} contains a version of Theorem \ref{CSTthm1.3} for the so-called $\Phi$-Laplacian. A more general problem patterned after
\begin{equation}\label{measprob}
-\Delta_p u=\mu\, u^{-\gamma}\;\;\mbox{in}\;\;\Omega,\;\; u>0\;\; \mbox{in}\;\;\Omega,\;\; u=0\;\;\mbox{on}\;\;\partial\Omega,   
\end{equation}
where $\mu$ denotes a non-negative bounded Radon measure on $\Omega$ while $\gamma\geq 0$, is thoroughly studied in \cite{DDO}; see also  \cite{DeCGOP} and the references therein.

Finally, as regards Problem \eqref{prob1} again, the papers \cite{KP,PS15,PS16,PW1} do not require Ambrosetti-Rabinowitz's condition $({\rm a}_5)$, while \cite{FS} establishes the existence of at least three weak solutions. Moreover, a possibly non-homogeneous elliptic operator is considered in \cite{PS16}, but $\lambda=1$.

The nice paper \cite{GST} investigates the problem
\begin{equation}\label{prob3}
\left\{
\begin{alignedat}{2}
-\Delta_p u & =\lambda u^{-\gamma}+u^{q-1}\;\;&&\mbox{in}\;\;\Omega, \\
u & >0 &&\mbox{in} \; \;\Omega, \\
u &= 0 &&\mbox{on} \; \;\partial\Omega,
\end{alignedat}
\right.
\end{equation}
where $0<\gamma<1$ and $1<p<q<p^*$. It should be noted that, here, contrary to above, \textit{the parameter $\lambda$ multiplies the singular term}. Combining known variational methods with a $C^{1,\alpha}(\overline{\Omega})$-regularity result \cite[Theorem 2.2]{GST} for solutions to \eqref{prob3} and a strong comparison principle \cite[Theorem 2.3]{GST}, the authors obtain the following
\begin{thm}[\cite{GST}, Theorem 2.1]\label{thm2}
Suppose $0<\gamma<1$ and $1<p<q<p^*$. Then there is $\Lambda>0$ such that \eqref{prob3} has:
\begin{itemize}
\item at least two ordered solutions in $C^1(\overline{\Omega})$ for every $\lambda\in\,]0,\Lambda[$,
\item at least one solution in $C^1(\overline{\Omega})$ when $\lambda=\Lambda$, and
\item no solutions once $\lambda>\Lambda$.
\end{itemize}
\end{thm}
The case $q=p^*$ is also studied and it is shown that $\gamma<1$ is a reasonable sufficient (and likely optimal) condition to get $C^1(\overline{\Omega})$-solutions of \eqref{prob3}.

If $p=2$ and, roughly speaking, $a\equiv-1$ while $f$ does not depend on $u$ then Problem \eqref{prob1} was fruitfully studied in \cite{DMO}.

We end this section by pointing out two very recent works, namely \cite{CGP}, which deals with possibly non-monotone singular reactions (see also \cite{GRS,H}, essentially based on sub-super-solution methods) and \cite{PW2}, devoted to singular equations driven by the $(p,q)$-Laplace operator $u\mapsto\Delta_p u+\Delta_q u$.
\subsection{Uniqueness}
Surprisingly enough, if $p\neq 2$, uniqueness of solutions looks a difficult matter, even for the model problem
\begin{equation}\label{prob4}
\left\{
\begin{alignedat}{2}
-\Delta_p u & =a(x)u^{-\gamma}\;\;&&\mbox{in}\;\;\Omega, \\
u & >0 &&\mbox{in} \; \;\Omega, \\
u &= 0 &&\mbox{on} \; \;\partial\Omega.
\end{alignedat}
\right.
\end{equation}
As observed in \cite{CST}, this is mainly caused by the fact that, in general, solutions do not belong to $W^{1,p}_0(\Omega)$ once $\gamma\geq 1$. The paper \cite{CST} provides two different results. The first one (Theorem 1.4) holds in star-shaped domains, while the other is the following
\begin{thm}[\cite{CST}, Theorem 1.5]
Assume that either $\gamma\leq 1$ and $a\in L^1(\Omega)$ or $\gamma>1$ and
\begin{itemize}
\item $a\in L^m(\Omega)$ for some $m>\frac{N}{p}$ if $1<p<N$, 
\item $a\in L^m(\Omega)$ with $m>1$ when $p=N$, and
\item $a\in L^1(\Omega)$ if $p>N$.
\end{itemize}
Then \eqref{prob4} possesses a unique distributional solution.
\end{thm}
We next point out that, for $\gamma\leq 1$, Theorem 3.4 of \cite{DDO} establishes the uniqueness of \textit{renormalized} solutions to \eqref{measprob}.

The situation becomes quite clear when $p=2$ and one seeks sufficiently regular solutions. Denote by $\varphi_1$ a positive eigenfunction corresponding to the first eigenvalue $\lambda_1$ of the problem $-\Delta u=\lambda u$ in $\Omega$, $u=0$ on $\partial\Omega$.
\begin{thm}[\cite{LM}, Theorems 1--2]
Let $p=2$ and let $a\in C^{0,\alpha}(\overline{\Omega})$ be positive. Then \eqref{prob4} has a unique solution $u\in C^{2,\alpha}(\Omega)\cap C^0(\overline{\Omega})$. Moreover,
\begin{itemize}
\item there exist $c_1,c_2>0$ such that $c_1\varphi_1^{2/(1+\gamma)}\leq u\leq c_2\varphi_1^{2/(1+\gamma)}$ in $\overline{\Omega}$,
\item $u\in H^1_0(\Omega)\iff\gamma<3$, and
\item $\gamma>1\implies u\not\in C^1(\overline{\Omega})$.
\end{itemize}
\end{thm}
See also the nice paper \cite{BGH}. As regards weak solutions, one has
\begin{thm}[\cite{CS}, Theorem 3.1]
Suppose $p=2$ and $a\in L^1(\Omega)$. Then \eqref{prob4} admits at most one solution belonging to $H^1_0(\Omega)$.
\end{thm}
Another uniqueness case occurs when $\gamma>1$.
\begin{thm}[\cite{CS}, Theorem 1.3]
If $p=2$, $\gamma>1$, and $a\in L^1(\Omega)$ then \eqref{prob4} possesses at most one solution $u\in H^1_{loc}(\Omega)$ such that $u^{(\gamma+1)/2}\in H^1_0(\Omega)$.
\end{thm}
\subsection{Equations with convective terms}
Consider the problem 
\begin{equation}\label{prob5}
\left\{
\begin{alignedat}{2}
-\Delta_p u & =f(x,u,\nabla u)+g(x,u)\;\;&&\mbox{in}\;\;\Omega,\\
u & >0 &&\mbox{in} \;\;\Omega,\\
u &= 0 &&\mbox{on} \;\;\partial\Omega,
\end{alignedat}
\right.
\end{equation}
where $p<N$ while $f:\Omega\times\R^+_0\times\R^N\to\R^+_0$ and $g:\Omega\times\R^+\to\R^+_0$ satisfy Carath\'eodory's conditions. In 2019, Liu, Motreanu, and Zeng established the existence of solutions $u\in W^{1,p}_0(\Omega)$ to \eqref{prob5} under the hypotheses below, where $\lambda_1$ stands for the first eigenvalue of $-\Delta_p$ in $W^{1,p}_0(\Omega)$.
\begin{itemize}
\item [$({\rm h}_1)$] There exist $c_0,c_1,c_2>0$ such that $c_1+c_2\lambda_1^{1-1/p}<\lambda_1$ and
$$f(x,t,\xi)\leq c_0+c_1 t^{p-1}+c_2|\xi|^{p-1}\;\;\forall\, (x,t,\xi)\in\Omega\times\R^+_0\times\R^N.$$
\item[$({\rm h}_2)$] $g(x,\cdot)$ is non-increasing on $(0,1]$ for all $x\in\Omega$ and $ g(\cdot,1)\not\equiv 0$.
\item[$({\rm h}_3)$] With appropriate $\theta\in\rm{int} (C^1_0(\overline{\Omega})_+)$, $\hat q>\max\{N,p'\}$, and $\epsilon_0 > 0$, the map $x\mapsto g(x,\epsilon \theta(x))$ belongs to $L^{\hat q}(\Omega) $ for any $\epsilon \in (0,\epsilon_0)$. 
\end{itemize}
Condition $({\rm h}_3)$ was previously introduced by Faraci and Puglisi \cite{FP}. It represents a natural generalization of $({\rm a}_1)$ in Section \ref{S2.1}.
\begin{thm}[\cite{LMZ}, Theorem 25]
Let $({\rm h}_1)$--$({\rm h}_3)$ be satisfied. Then \eqref{prob5} has a solution $u\in\rm{int}(C^1_0(\overline{\Omega})_+)$.
\end{thm}
We think worthwhile to sketch the main ideas of the proof. For every fixed $w\in C^1_0(\overline{\Omega})$, an intermediate problem, where $\nabla w$ replaces $\nabla u$ in $f(x,u,\nabla u)$ and the singular term remains unchanged, is considered. The authors construct a positive sub-solution $\underline{u}\in\rm{int}(C^1_0(\overline{\Omega})_+)$ independently of $w$ and show the existence of a solution greater than $\underline{u}$. If ${\cal S}(w)$ denotes the set of such solutions then, via suitable properties of the multi-function $w\mapsto{\cal S}(w)$, it is proved that the map $\Gamma$, which assigns to each $w$ the minimal element of ${\cal S}(w)$, is completely continuous. Now, Leray–Schauder's alternative principle applied to $\Gamma$ yields a solution $u\in\rm{int}(C^1_0(\overline{\Omega})_+)$ to \eqref{prob5}.

The recent paper \cite{GMM}, partially patterned after \cite{LMZ}, treats the Robin problem
\begin{equation}\label{prob5bis}
\left\{
\begin{alignedat}{2}
-{\rm div}\,A(\nabla u) & =f(x,u,\nabla u)+g(x,u)\;\;&&\mbox{in}\;\;\Omega,\\
u & >0 &&\mbox{in} \;\;\Omega,\\
\frac{\partial u}{\partial\nu_A}+\beta u^{p-1} & =0 &&\mbox{on} \;\;\partial\Omega,
\end{alignedat}
\right.
\end{equation}
where $A:\R^N\to\R^N$ denotes a continuous strictly monotone map having suitable properties, which basically stem from Lieberman's nonlinear regularity theory \cite{L} and Pucci-Serrin's maximum principle \cite{PuSe}. By the way, the conditions on $A$ include classical non-homogeneous operators as, e.g., the $(p,q)$-Laplacian. Moreover,
$\beta$ is a positive constant while $\frac{\partial}{\partial\nu_A}$ indicates the co-normal derivative associated with $A$.  If $p=2$ then a uniqueness result is also presented; cf. \cite[Theorem 4.2]{GMM}.

The special case $A(\xi):=|\xi|^{p-2}\xi$, $g(x,t):=t^{-\gamma}$ for some $0<\gamma<1$, and $\beta=0$ (which reduces \eqref{prob5bis} to a Neumann problem) has been investigated in \cite{PRR20} without imposing any global growth condition on $t\mapsto f(x,t,\xi)$. Instead, a kind of oscillatory behavior near zero is taken on. For such an $f$, the work \cite{PaZ} establishes the existence of a solution $u\in C^1_0(\overline{\Omega})$ to the parametric problem
\begin{equation*}\label{prob5ter}
\left\{
\begin{alignedat}{2}
-{\rm div}\,A(\nabla u) & =f(x,u,\nabla u)+\lambda u^{-\gamma}\;\; && \mbox{in}\;\;\Omega,\\
u & >0 &&\mbox{in} \;\;\Omega,\\
u & =0 &&\mbox{on} \;\;\partial\Omega,
\end{alignedat}
\right.
\end{equation*}
provided $\lambda>0$ is small enough.

Finally, the very recent paper \cite{CGSS} treats $\Phi$-Laplacian equations with strongly singular reactions perturbed by gradient terms.
\section{Problems on the whole space}\label{whole}
\subsection{The case $p=2$}
Let $N\geq 3$, let $a:\R^N\to\R^+_0$ be nontrivial measurable, and let $\gamma>0$. The simplest singular elliptic problem in the whole space writes as
\begin{equation}\label{probwholespace}
\left\{
\begin{alignedat}{2}
-\Delta u & =a(x)u^{-\gamma}\;\;&&\mbox{in}\;\;\R^N,\\
u & >0 &&\mbox{in} \;\;\R^N.\\
\end{alignedat}
\right.
\end{equation}
Sometimes it is also required that $u(x)\to 0$ as $|x|\to\infty$. Since the pioneering papers \cite{KS,Dal,E,LS1}, some existence and uniqueness results concerning \eqref{probwholespace} have been published. We refer the reader to the monograph \cite{GR} for a deep account. Roughly speaking, four basic questions can be identified:
\begin{itemize}
\item find the right hypotheses on $a$. Usually, $a\in C^{0,\alpha}_{loc} (\R^N)_+$ as well as
$$\int_1^\infty r\max_{|x|=r}a(x)\,{\rm d}r<\infty$$
(cf. condition $({\rm a}_8)$ below) guarantee both existence and uniqueness of solutions $u\in C^{2,\alpha}_{loc} (\R^N)$.
\item replace $u^{-\gamma}$ with a function $f(u)$ such that $\lim_{t\to 0^+}f(t)=\infty$. This was done in \cite{LS2,Z} for decreasing $f$. Later on, also non-monotone singular reactions were fruitfully treated \cite{CR,GS2,GMS}.
\item put convective terms on the right-hand side. For equations driven by the Laplacian, a good reference is  \cite[Section 9.8]{GR}; cf. in addition \cite{GR1,GSi}. 
\item generalize the left-hand side of the equation. The case of a second-order uniformly elliptic operator is treated in \cite{ChK,Chab}, while \cite{AGM} deals with $u\mapsto-\Delta u+c(x)u$, where $c\in L^\infty_{loc}(\R^N)_+$.
\end{itemize}
The equation of Problem \eqref{probwholespace} arises in the boundary-layer theory  of viscous fluids \cite{CF,CN2,CN3} and is called \textit{Lane-Emden-Fowler equation}. Its importance in scientific applications has by now been widely recognized; see, e.g., \cite{Fo}.
\subsection{Existence and multiplicity}
To the best of our knowledge, the first paper treating singular $p$-Laplacian equations on the whole space is that of Goncalves and Santos \cite{GS1}, published in 2004. The authors consider the problem
\begin{equation}\label{prob6}
\left\{
\begin{alignedat}{2}
-\Delta_p u & =a(x)f(u)\;\;&&\mbox{in}\;\;\R^N,\\
u & >0 &&\mbox{in} \;\;\R^N,\\
u(x) & \to 0 &&\mbox{as}\;|x|\to\infty, 
\end{alignedat}
\right.
\end{equation}
where $a\in C^0(\R^N)_+$ is radially symmetric while $f\in C^1(\R^+,\R^+)$, and assume that:
\begin{itemize}
\item[$({\rm a}_6)$] the function $t\mapsto\frac{f(t)}{t^{p-1}}$ is non-increasing on $\R^+$.
\item[$({\rm a}_7)$] $\displaystyle{\liminf_{t\to 0^+}}f(t)>0$ as well as $\displaystyle{\lim_{t\to\infty}}\,\frac{f(t)}{t^{p-1}}=0$.
\item[$({\rm a}_8)$] if $\Phi(r):=\max_{|x|=r}a(x)$, $r>0$, then
\begin{equation*}
\begin{alignedat}{2}
& 0<\int_1^\infty [r\Phi(r)]^{\frac{1}{p-1}}{\rm d}r<\infty && \mbox{for}\; 1<p\leq 2,\\
& 0<\int_1^\infty r^{\frac{(p-2)N+1}{p-1}}\Phi(r)\,{\rm d}r<\infty\;\; && \mbox{for}\; p>2.
\end{alignedat}
\end{equation*}
\end{itemize}
\begin{thm}[\cite{GS1}, Theorem 1.1]\label{thmGS1}
Under $({\rm a}_6)$--$({\rm a}_8)$, Problem \eqref{prob6} admits:
\begin{itemize}
\item a radially symmetric solution $u\in C^1(\R^N)\cap C^2(\R^N\setminus\{0\})$ when $p<N$.
\item no radially symmetric solution in $C^1(\R^N)\cap C^2(\R^N\setminus\{0\})$ if $p\geq N$.
\end{itemize}
\end{thm}
The proof exploits fixed point arguments, the shooting method, and sub-super-solution techniques.

One year later, Covei \cite{Co1} did not assume $a$ radially symmetric but locally H\"{o}lder continuous and positive, replaced conditions $({\rm a}_6)$--$({\rm a}_7)$ with those below, and obtained similar results. See also \cite{Co2}, where the asymptotic behavior of solutions is described.
\begin{itemize}
\item[$({\rm a}_6')$] The function $t\mapsto\frac{f(t)}{(t+\beta)^{p-1}}$ turns out decreasing on $\R^+$ for some $\beta>0$.
\item[$({\rm a}_7')$] $\displaystyle{\lim_{t\to 0^+}} \frac{f(t)}{t^{p-1}}=\infty$ and $f(t)\leq c$ for any $t$ large enough.
\end{itemize}
The work \cite{LGL} treats the parametric problem
\begin{equation}\label{prob7}
\left\{
\begin{alignedat}{2}
-\Delta_p u & =a(x)u^{-\gamma}+\lambda b(x)u^{q-1}\;\; && \mbox{in} \;\;\R^N,\\
u & >0 &&\mbox{in} \;\;\R^N,
\end{alignedat}
\right.
\end{equation}
where $1<p<N$, $0<\gamma<1$, $\lambda>0$, $\max\{p,2\}<q<p^*$, and the coefficients fulfill
\begin{equation}\label{hpab}
a\in L^{\frac{p^*}{p^*-(1-\gamma)}}(\R^N)_+,\;\;a\not\equiv 0,\;\;  b\in L^{\frac{p^*}{p^*-q}}(\R^N), \;\; b>0.  
\end{equation}
\begin{thm}[\cite{LGL}, Theorem 1.2]
If \eqref{hpab} holds then there exists $\Lambda>0$ such that \eqref{prob7} possesses
\begin{itemize}
\item at least two solutions in ${\cal D}^{1,p}_0(\R^N)$ for every $\lambda\in\, ]0,\Lambda[$,
\item at least one solution belonging to ${\cal D}^{1,p}_0(\R^N)$ when $\lambda=\Lambda$, and
\item no solutions once $\lambda>\Lambda$.
\end{itemize}
\end{thm}
It may be of interest to point out that this result is proved by combining sub-super-solution methods with the mountain pass theorem for continuous functionals.
\begin{rmk}
If $b\equiv 0$ then Problem \eqref{prob7} reduces to a well-known one, very important in scientific applications; cf. \cite[Remark 2.2]{CP}.
\end{rmk}
A meaningful case occurs when $a,b:\R^N\to\R^+_0$ turn out nonzero locally H\"{o}lder continuous functions. In fact, define 
\begin{equation}\label{defM}
M(x):=\max\{a(x),b(x)\},\quad x\in\R^N.    
\end{equation}
From \cite[Remarks 1--2]{Sa} it follows
\begin{lemma}\label{aux1}
Suppose that $p<N$, the functions $a,b:\R^N\to\R^+_0$ are nontrivial and locally H\"{o}lder continuous, while $({\rm a}_8)$ holds with $M$ in place of $a$. Then the problem
\begin{equation}
\left\{
\begin{alignedat}{2}
-\Delta_p w & =M(x)\;\;&&\mbox{in}\;\;\R^N,\\
w & >0 &&\mbox{in}\;\;\R^N,\\
w(x) & \to 0 &&\mbox{as}\;|x|\to\infty
\end{alignedat}
\right.
\end{equation}
admits a solution $w_M\in C^{1,\alpha}_{loc}(\R^N)$ for suitable $\alpha\in\,]0,1[$.
\end{lemma}
Via sub-super-solution techniques, Lemma \ref{aux1} gives rise to
\begin{thm}[\cite{Sa}, Theorem 1.1]
Let $\gamma>0$, let $p<q$, and let $M$ be given by \eqref{defM}. Under the assumptions of Lemma \ref{aux1}, there exists $\lambda^*>0$ such that \eqref{prob7} has:
\begin{itemize}
\item at least one solution $u\in C^1(\R^N)$ for every $0\leq\lambda<\lambda^*$. Moreover, $u(x)\to 0$ as $|x|\to\infty$.
\item no solution once $\lambda>\lambda^*$.
\end{itemize}
\end{thm}
This result was next generalized under various aspects by the same author and Rezende \cite{RS}; cf. also \cite{CP}.

Finally, infinite semi-positone problems, i.e., $\lim_{t\to 0^+} f(t)=-\infty$, were fruitfully investigated in \cite{DS}. Precisely, given $a\in L^\infty(\R^N)$ and $f\in C^0(\R^+)$, consider the problem
\begin{equation}\label{prob8}
\left\{
\begin{alignedat}{2}
-\Delta_p u & =\lambda a(x)f(u)\;\;&&\mbox{in}\;\;\R^N,\\
u & >0 &&\mbox{in} \;\;\R^N,\\
u(x) & \to 0 &&\mbox{as}\;|x|\to\infty, 
\end{alignedat}
\right.
\end{equation}
where $\lambda>0$, $1<p<N$. The following conditions will be posited.
\begin{itemize}
\item[$({\rm a}_9)$] There exists $\gamma\in\,]0,1[$ such that $\displaystyle{ \lim_{t\to 0^+}}t^{\gamma}f(t)=c_0\in\R^-$.
\item[$({\rm a}_{10})$] $\displaystyle{\lim_{t\to\infty}f(t)=\infty}$ but $\displaystyle{\lim_{t\to\infty}}\,\frac{f(t)}{t^{p-1}}=0$.
\item[$({\rm a}_{11})$] $\displaystyle{\inf_{|x|=r}}a(x)>0$ for all $r>0$ and $0<a(x)<\frac{C_0}{|x|^\sigma}$ in $\R^N\setminus\{0\}$ with suitable $C_0>0$, $\sigma>N+\gamma\frac{N-p}{p-1}$.
\end{itemize}
\begin{thm}[\cite{DS}, Theorem 1.4]
If $({\rm a}_9)$--$({\rm a}_{11})$ hold and $\lambda$ is sufficiently large then \eqref{prob8} has a solution in $C^{1,\alpha}_{loc}(\R^N)$. 
\end{thm}
\subsection{Uniqueness}
As far as we know, uniqueness has been addressed only in \cite[Remark 1.2]{GS1} and \cite[Section 2]{Co1} under the key assumption $({\rm a}_6')$ above. The arguments of both papers rely on a famous result by Diaz and Saa \cite{DSa}. Theorem 1.3 of \cite{CDS} contains a nice idea to achieve uniqueness for singular problems in exterior domains.
\subsection{Equations with convective terms}
To the best of our knowledge, there is only one paper concerning singular quasi-linear elliptic equations in the whole space and with convective terms, namely \cite{GG}. It treats the problem
\begin{equation}\label{GGprob}
\left\{
\begin{alignedat}{2}
-{\rm div} \, A(\nabla u) &= f(x,u) + g(x,\nabla u) \quad &&\mbox{in} \;\; \R^N, \\
u &> 0 \quad &&\mbox{in} \;\; \R^N,
\end{alignedat}
\right.
\end{equation}
where $N \geq 2$ and $1<p<N$. The differential operator $u\mapsto{\rm div}\, A(\nabla u)$ is as in \eqref{prob5bis}, while $f:\R^N\times\R^+ \to\R^+_0$ and $g:\R^N\times\R^N\to\R^+_0$ fulfill Carathéodory's conditions. Moreover,
\begin{equation}\label{hypf}
\begin{split}
\liminf_{t\to 0^+}f(x,t)>0\;\;\mbox{uniformly with respect to}\;\; x\in 
B_\sigma(x_0), \\
f(x,t)\leq h(x)t^{-\gamma}\;\mbox{in}\;\R^N\times\R^+,\;\mbox{where}\; h\in L^1(\R^N)\cap L^\eta(\R^N),
\end{split}
\end{equation}
and
\begin{equation}\label{hypg}
g(x,\xi)\leq k(x)|\xi|^r\;\mbox{in}\;\R^N\times\R^N,\;\mbox{with}\; k\in L^1(\R^N)\cap L^\theta(\R^N).
\end{equation}
Here, $x_0\in\R^N$, $\sigma\in ]0,1[$, $\gamma\geq 1$, $r\in [0,p-1[$, as well as
\begin{equation}\label{hyphkr}
\eta>(p^*)',\quad\theta > \left(\frac{1}{(p^*)'}-\frac{r}{p}\right)^{-1}.
\end{equation}
\begin{thm}[\cite{GG}, Theorem 1.2]\label{GGthm}
Under \eqref{hypf}--\eqref{hyphkr}, there exists a distributional solution $u\in W^{1,p}_{loc}(\R^N)$ to \eqref{GGprob} such that
$\essinf_K u > 0$ for every compact set $K\subseteq\R^N$.
\end{thm}
To prove this result, the authors first solve some auxiliary problems, obtained by shifting the singular term and working in balls, via  sub-super-solution techniques. A compactness result, jointly with a fine local energy estimate on super-level sets of solutions, then yields the conclusion.
\small{
\subsection*{Acknowledgments}
U. Guarnotta and S.A. Marano were supported by the research project `MO.S.A.I.C.' PRA 2020--2022 `PIACERI' Linea 2 (S.A. Marano) and Linea 3 (U. Guarnotta)  of the University of Catania.\\
U. Guarnotta, R. Livrea, and S.A. Marano were supported by the research project PRIN 2017 `Nonlinear Differential Problems via Variational, Topological and Set-valued Methods' (Grant No. 2017AYM8XW) of MIUR.\\
U. Guarnotta also acknowledges the support of the GNAMPA-INdAM Project CUP\_E55F22000270001.

\noindent
\textbf{Conflict of interest statement:} the authors state no conflict of interest.
}

\end{document}